\documentclass[11pt,twoside]{amsart} 
\usepackage[koi8-r]{inputenc}
\usepackage[english]{babel}
\usepackage{amssymb}
\usepackage{amsmath}
\advance \topmargin by -\headheight
\advance \topmargin by -\headsep
\evensidemargin \oddsidemargin
\date{}
\newtheorem{conjecture}{Conjecture}[section]
\newtheorem{theorem}[conjecture]{Theorem}
\newtheorem{conj-corollary}{``Corollary''}[conjecture]
\newtheorem{corollary}[conjecture]{Corollary}
\newtheorem{lemma}[conjecture]{Lemma}
\newtheorem{proposition}[conjecture]{Proposition}
\title{Representations of field automorphism groups}
\author{M.Rovinsky} 
\begin{document} 
\begin{abstract}
This is a common introduction to {\tt math.AG/0011176, math.RT/0101170, 
math.RT/0306333, math.RT/0506043, math.RT/0601028}. 

In these papers one studies the automorphism group $G$ of an extension $F/k$ 
of algebraically closed fields, especially in the case of countable 
transcendence degree and characteristic zero, its smooth linear 
and semi-linear representations, and their relations to algebraic 
geometry (birational geometry, motives, differential forms and sheaves). 

Compared to the above references there are some new results 
including \begin{itemize} \item a description of a separable closure 
of an extension of transcendence degree one of an algebraically 
closed field (Proposition \ref{add-mult}, p. \pageref{add-mult});
\item a ``K\"{u}nneth formula'' for the products with curves;
\item the semi-simplicity of the $G$-module $\Omega^n_{F/k,{\rm reg}}$ 
of regular differential forms of top degree. \end{itemize} 
\end{abstract}
\maketitle 

The study of field automorphism groups is an old subject. Without 
any attempt of describing its complicated history, let me just 
mention that many topological groups are field automorphism 
groups. Besides the usual Galois groups we meet here (discrete, 
$p$-adic for $p<\infty$, or finite adelic) groups of points of 
algebraic groups. 

Let $F/k$ be a field extension of countable (this will be the principal 
case) or finite transcendence degree $n$, $0\le n\le\infty$, and 
$G=G_{F/k}$ be its automorphism group. Following \cite{jac,pss,sh,ih} 
(and generalizing the case \cite{krull-gal} of algebraic extension), 
consider $G$ as a topological group with the base of open subgroups 
given by the stabilizers of finite subsets of $F$. Then $G$ is a 
totally disconnected Hausdorff group, and for any intermediate 
subfield $k\subseteq K\subseteq F$ the topology on $G_{F/K}$ 
coincides with the restriction of the topology on $G$. There are maps: 
from the set of intermediate subfields in $F/k$ to the set of closed 
subgroups of $G$, $K\mapsto G_{F/K}:={\rm Aut}(F/K)$, and from the set 
of closed subgroups in $G$ to the set of intermediate subfields in 
$F/k$, $H\mapsto F^H$. They are mutually inverse to each other in the 
Galois extension case. If $n<\infty$ then $G$ is locally compact. 

Following the very general idea, not only in Mathematics, that a 
``sufficiently symmetric'' system is determined by a representation 
of its symmetry group, one tries to compare various ``geometric 
categories over $k$'' with various categories of representations 
of $G$, 

To ensure that the representation theory of $G$ is rich enough, $F$ 
should be ``big enough'', e.g. algebraically closed. So $F$ is ``the 
function field of the universal tower of $n$-dimensional $k$-varieties'', 
if $n<\infty$. In that case each perfect subfield $L$ of $F$ containing 
$k$ is the fixed field of the subgroup $G_{F/L}$ of $G$, and $G$ 
contains, in particular, the groups $G_{L/k}$ as its sub-quotients. 

Usually (unless it is not stated otherwise), $k$ will 
be algebraically closed (in order to avoid already 
complicated enough Galois theory) of characteristic zero. 

One of the main motivations is the calculation of integrals of 
meromorphic differential forms $\omega$ on projective complex 
varieties. To calculate such an integral, one can transfer 
$\omega$ to other varieties via correspondences. In coordinates 
this looks as an algebraic change of variables. We may suppose that 
all function fields are contained in a common field $F$. Then 
the problem of description of the properties of the (iterated) 
integrals of $\omega$ (of $\omega_1,\dots\omega_N$) becomes related 
with determining the structure of the $G$-submodule in the algebra 
of K\"{a}hler differentials $\Omega^{\bullet}_{F/k}$ (resp., in 
$\Omega^{\bullet}_{F/k}\otimes_k\dots\otimes_k\Omega^{\bullet}
_{F/k}$) generated by $\omega$ (resp., by $\omega_1\otimes\dots
\otimes\omega_N$). 

For example, the irreducible subquotients of the $G$-module 
$\Omega^1_{F/k,{\rm closed}}$ of closed K\"{a}hler differentials 
are related to the simple algebraic commutative groups 
over $k$, cf. Proposition \ref{1-forms}. 

In the opposite direction, to clarify relations, so far conjectural, 
between the motives and the cohomologies, one has to link the most 
interesting (from the geometric point of view) representations -- 
admissible and ``homotopy invariant'' -- and the K\"{a}hler differentials. 
Conjecturally, the irreducible ones among them are contained in the 
algebra of differential forms $\Omega^{\bullet}_{F/k}$, if $n=\infty$. 

\subsection{Some general notations, conventions and goals} 
\label{some-notation}
Let $F/k$ be an extension of countable or finite transcendence 
degree $n$, $1\le n\le\infty$, of algebraically closed fields of 
characteristic zero (by default), and $G=G_{F/k}$ be its 
automorphism group endowed with the above topology. 

We study the structure of $G$, its linear and semi-linear representations 
(with open stabilizers), and their relations to algebraic geometry 
(birational geometry, motives, differential forms and sheaves) and 
to automorphic representations. 
In particular, we look for analogues of known results for 
$p$-adic (and more generally, locally compact) groups. 

\section{Structure of $G$} \label{structure-of-G}
It is well-known (\cite{jac,pss,sh,ih}) that the group $G$ 
is locally compact if and only if $n<\infty$. 
\begin{theorem}[\cite{repr}] \label{top-simpl} \begin{enumerate} 
\item The subgroup $G^{\circ}$ of $G$, generated by the compact 
subgroups, is open and topologically simple, if $n<\infty$. 
If $n=\infty$ then $G^{\circ}$ is dense in $G$. 
\item \label{prosto-besk} Any closed normal proper subgroup of 
$G$ is trivial, if $n=\infty$, i.e. $G$ is topologically simple. 
\end{enumerate} \end{theorem}

{\sc Remarks.} 1. In fact, Theorem \ref{top-simpl} holds for 
any extension $F/k$ of algebraically closed fields of arbitrary 
characteristic, cf. \cite{repr}. Moreover, if $n=1$ and ${\rm char}(k)
\neq 0$ then the separable closure of $k(x)$ in $F$ is generated by the 
$G^{\circ}$-orbit of $x$ for any $x\in F\backslash k$, cf. \cite{obz}. 

2. An argument of \cite{lasc} shows that $G$ is simple as a discrete 
group provided that transcendence degree $F$ over $k$ is not countable. 

3. If $n<\infty$, the left $G$-action on the one-dimensional 
oriented ${\mathbb Q}$-vector space of right-invariant measures 
on $G$ gives rise to a surjective homomorphism, the modulus, 
$\chi:G\longrightarrow{\mathbb Q}^{\times}_+$, which is 
trivial on $G^{\circ}$. However, I do not know even, whether 
the discrete group $\ker\chi/G^{\circ}$ is trivial. If it 
is trivial for $n=1$ then it is trivial in general, cf. \cite{repr}. 

\subsection{Closed, open and maximal proper subgroups; Galois theories}
The classical morphism $\beta:\{\mbox{subfields $F$ over $k$}\}
\hookrightarrow\{\mbox{closed subgroups of $G$}\}$, given by 
$K\mapsto G_{F/K}$, is injective, inverts the inclusions, 
transforms the compositum of subfields to the intersection of 
subgroups, and respects the units: $k\mapsto G$. The image of 
$\beta$ is stable under the passages to sup-/sub- groups with 
compact quotients; $\beta$ identifies the subfields over 
which $F$ is algebraic with the compact subgroups of $G$ 
(\cite{jac,pss,sh,ih}). 

In particular, the proper subgroups in 
the image of $\beta$ are the compact subgroups if $n=1$. 

The map $H\mapsto F^H$, left inverse of $\beta$, inverts the order, 
but does not respect the monoid structure. 

In \cite{max}, a morphism of partially ordered commutative associative 
unitary monoids (transforming the intersection of subgroups to 
the algebraic closure of the compositum of subfields) 
$$\alpha:\{\mbox{open subgroups of $G$}\}\longrightarrow\!\!\!\!\rightarrow
\left\{\begin{array}{c}\mbox{algebraically closed subfields of $F$}\\ 
\mbox{of finite transcendence degree over $k$}\end{array}\right\}$$ 
is constructed, It is determined uniquely by the condition $G_{F/\alpha(U)}
\subseteq U$ and the transcendence degree of $\alpha(U)$ over $k$ is minimal. 

\vspace{4mm}

It is shown in \cite{max} that for any non-trivial algebraically closed 
extension $L\neq F$ of $k$ of finite transcendence degree in $F$ the 
normalizer in $G$ of $G_{F/L}$ (which is evidently open) is maximal 
among the proper subgroups of $G$. In the case $n=\infty$ any maximal 
open proper subgroup is of this type. 

As a consequence, one gets a complete, though not very explicit, 
Galois theory of algebraically closed extensions of countable 
transcendence degree (a question of Krull, \cite{krull}), i.e. 
a construction of all subgroups $H$ of $G$ coincident with the 
automorphism groups of $F$ over the fixed subfields $F^H$. 

Another type of closed non-open maximal proper subgroups is given by 
the stabilizers of rank one discrete valuations in the case of arbitrary 
transcendence degree. They are useful in relating representations of $G$ 
to functors on categories of smooth $k$-varieties, cf. \S\ref{usl-glob-porozhd}. 

\subsection{Automorphisms of $G$}
The group $G$ is quite rigid in the sense that the group of its continuous 
automorphisms is ``of the same size'' as $G$. Namely, it coincides with 
the group of field automorphisms of $F$ preserving $k$. If $n\ge 2$ this 
follows from results of F.A.Bogomolov. If $n=1$ this is shown in 
\cite{iso-abs}. 

It would be highly interesting to identify the class of ``rational'' 
representations of $G$, i.e. those whose isomorphism class does not 
change under any continuous automorphism of $G$. In particular, if 
$L\subset F$ is a field of automorphic functions (of all levels) 
the functor $H^0(G_{F/L},-)$ should relate representations of $G$ 
to automorphic representations. 
 
\section{How to translate geometric questions to the language of 
representation theory?}
Depending on type of geometric questions we shall consider 
one of the following threer categories of representations of $G$: 
${\mathcal S}m_G\supset{\mathcal I}_G\supset{\mathcal A}dm$, 
roughly corresponding to birational geometry over $k$, birational 
motivic questions (like on the structure of Chow groups of 0-cycles) 
and ``finite-dimensional'' birational motivic questions (such as 
description of ``classical'' motivic categories). 

\subsection{${\mathcal S}m_G$} Usually an ``algebro-geometric datum'' 
$D$ over $F$, a universal domain over $k$ in the sense of Weil, 
consists of a finite number of polynomial equations involving a 
finite number of coefficients $a_1,\dots,a_N\in F$, and the group 
$G$ acts on the set of ``similar'' data. Then the stabilizer of 
$D$ in $G$ is open, since contains $G_{F/k(a_1,\dots,a_N)}$. 

For a $k$-variety $X$, its $F$-subvarieties are examples of such data. 

In particular, the ${\mathbb Q}$-vector space ${\mathbb Q}[X(F)]$ 
of 0-cycles on $X\times_kF$ is a $G$-module, 
Such representation is huge, but this is just a starting point. 

Note that it is smooth, i.e. its stabilizers are open, so all 
representations we are going to consider will be smooth. 

Conversely, as it follows from \cite{repr}, Lemma 3.3, any smooth 
representation of $G$ with cyclic vector is a quotient of the $G$-module 
${\mathbb Q}[\{k(X)\stackrel{/k}{\hookrightarrow}F\}]$ of ``generic'' 
0-cycles on $X_F$ (equivalently, formal ${\mathbb Q}$-linear combinations 
of embeddings of the function field $k(X)$ into $F$ over $k$), i.e., 
0-cycles outside of the union of the divisors on $X$ defined over $k$, 
for an appropriate irreducible variety $X$ of dimension $\le n$ over $k$. 

{\sc Remarks.} 1. One has ${\mathbb Q}[X(F)]= 
\bigoplus_{x\in X}{\mathbb Q}[\{k(x)\stackrel{/k}{\hookrightarrow}F\}]$, 
so ${\mathbb Q}[X(F)]$ reflects rather the class of $X$ in the Grothendieck 
group $K_0(Var_k)$ of partitions of varieties over $k$ than $X$ itself. 

2. It is not clear, whether the birational type of $X$ is determined 
by the $G$-module of generic 0-cycles on $X_F$. There exist pairs of 
non-birational varieties $X$ and $Y$, whose $G$-modules of generic 
0-cycles have the same irreducible subquotients, cf. \cite{repr}. 
Namely, $X=Z\times{\mathbb P}^1$ and $Y=Z'\times{\mathbb P}^1$, where 
$Z'$ is a twofold cover of $Z$, is such a pair. 
What is in common between $X$ and $Y$ in this example, is that their 
primitive motives (see below) coincide (and vanish). 

However, one can extract ``birational motivic'' invariants ``modulo 
isogenies'', such as ${\rm Alb}(X)$, $\Gamma(X,\Omega^{\bullet}_{X/k})$, 
out of ${\mathbb Q}[\{k(X)\stackrel{/k}{\hookrightarrow}F\}]$, cf. Theorem 
\ref{mnogo-o-I} (\ref{left-right-adj}--\ref{pse-mot-filt}). 

\vspace{4mm}

Denote by ${\mathcal S}m_G$ the category of 
smooth representations of $G$ over ${\mathbb Q}$. 

It follows from the topological simplicity of $G$ (Theorem 
\ref{top-simpl}) that in the case $n=\infty$ any finite-dimensional 
smooth representation of $G$ is trivial.

\subsection{${\mathcal A}dm$} \label{dopu-pred} Now consider 
a more concrete geometric category: the category of motives. 

(Effective) pure covariant motives are pairs $(X,\pi)$ 
consisting of a smooth projective variety $X$ over $k$ with 
irreducible components $X_j$ and a projector $\pi=\pi^2\in
\bigoplus_jB^{\dim X_j}(X_j\times_kX_j)$ in the algebra 
of correspondences on $X$ modulo numerical equivalence. 
The morphisms are defined by ${\rm Hom}((X',\pi'),(X,\pi))=
\bigoplus_{i,j}\pi_j\cdot B^{\dim X_j}(X_j\times_kX'_i)\cdot\pi'_i$. 
The category of pure covariant motives carries an additive 
and a tensor structures: $$(X',\pi')\bigoplus(X,\pi):=
(X'\coprod X,\pi'\oplus\pi),\qquad(X',\pi')\otimes(X,\pi)
:=(X'\times_kX,\pi'\times_k\pi).$$ 
A {\sl primitive $q$-motive} is a pair $(X,\pi)$ as above with 
$\dim X=q$ and ${\rm Hom}(Y\times{\mathbb P}^1,(X,\pi)):=\pi\cdot 
B^q(X\times_kY\times{\mathbb P}^1)=0$ for any smooth projective 
variety $Y$ over $k$ with $\dim Y<q$. For instance, due to 
the Lefschetz theorem on $(1,1)$-classes, the category of the 
pure primitive $1$-motives is equivalent to the category 
of abelian varieties over $k$ with morphisms tensored with 
${\mathbb Q}$. It is a result of Jannsen \cite{jan} that pure 
motives form a semi-simple abelian category, and it follows from 
this that any pure motive admits a ``primitive'' decomposition 
$\bigoplus_{i,j}M_{ij}\otimes{\mathbb L}^{\otimes i}$, 
where $M_{ij}$ is a primitive $j$-motive and ${\mathbb L}=
({\mathbb P}^1,{\mathbb P}^1\times\{0\})$ is the Lefschetz motive. 

\vspace{4mm}

{\sc Definition.} A representation $W$ of a topological group is called 
{\sl admissible} if it is smooth and the fixed subspaces $W^U$ are 
finite-dimensional for all open subgroups $U$. 

\vspace{5mm}

Denote by ${\mathcal A}dm$ the category of admissible representations 
of $G$ over ${\mathbb Q}$. 

\begin{theorem}[\cite{repr}] \label{adm-ab} ${\mathcal A}dm$ 
is a Serre subcategory in ${\mathcal S}m_G$. \end{theorem}
In other words, ${\mathcal A}dm$ is abelian, stable under taking 
sub{\sl quotients} (this is the point in the case $n=\infty$!) 
in the category of representations of $G$, and under taking 
extensions in ${\mathcal S}m_G$. 

\begin{theorem}[\cite{repr}] \label{vpolne-strog} There is a 
fully faithful functor ${\mathbb B}^{\bullet}$ when $n=\infty$: 
$$\left\{\begin{array}{c} \mbox{{\rm pure covariant motives over $k$}}
\end{array}\right\}
\stackrel{{\mathbb B}^{\bullet}}{\longrightarrow}
\left\{\begin{array}{c} \mbox{{\rm graded semi-simple admissible}}\\
\mbox{$G$-{\rm modules of finite length}}\end{array}\right\}.$$
The grading corresponds to powers of the motive 
${\mathbb L}$ in the ``primitive'' decomposition above. \end{theorem}
Roughly speaking, the functor ${\mathbb B}^{\bullet}$ is defined by 
spaces of 0-cycles defined over $F$ modulo ``numerical 
equivalence over $k$''. More precisely, ${\mathbb B}^{\bullet}=
\oplus^{{\rm graded}}_j\lim\limits_{_L\longrightarrow}{\rm Hom}
\left([L]^{{\rm prim}}\otimes{\mathbb L}^{\otimes j},-\right)$ is a  
graded direct sum of pro-representable functors. Here $L$ runs over all  
subfield of $F$ of finite type over $k$, and $[L]^{{\rm prim}}$ is the  
quotient of the motive of any smooth projective model of $L$ over $k$ by  
the sum of all submotives of type $M\otimes{\mathbb L}$ for all 
effective motives $M$. 

\vspace{4mm}

{\sc Examples.} The motive of the point ${\bf Spec}(k)$ is sent to the
trivial representation ${\mathbb Q}$ in degree 0. The motive of a smooth 
proper curve $C$ over $k$ is sent to ${\mathbb Q}\oplus J_C(F)/J_C(F)
\oplus{\mathbb Q}[1]$, where $J_C$ is the Jacobian of $C$ and 
${\mathbb Q}[1]$ denotes the trivial representation in degree 1. 

\vspace{4mm}

So, this inclusion is already a good reason to study admissible 
representations. 

\vspace{4mm}

Moreover, it is expected that 
\begin{conjecture} \label{conj-equi} The functor 
${\mathbb B}^{\bullet}$ is an equivalence of categories. \end{conjecture} 

Of course, it would be more interesting to describe in a similar way 
the abelian category $\mathcal{MM}$ of mixed motives over $k$, whose 
semi-simple objects are pure. This is one more reason to study 
the category ${\mathcal A}dm$ of admissible representations of $G$. 

\begin{proposition}[\cite{repr}] \label{ext-calc-ab} 
Assuming $n=\infty$, for any $W\in{\mathcal A}dm$, any abelian variety 
$A$ over $k$ and, conjecturally, for any effective motive $M$ one has \\
\begin{tabular}{ll} ${\rm Ext}^{>0}_{{\mathcal A}dm}({\mathbb Q},W)=0$ 
& ${\rm Ext}^{>0}_{\mathcal{MM}}({\mathbb Q},M)=0$ \\ 
${\rm Ext}^1_{{\mathcal A}dm}(\frac{A(F)}{A(k)},W)=
\frac{{\rm Hom}_{{\mathbb Z}}(A(k),W^G)}{{\rm Hom}_G(A(F)/A(k),W/W^G)}$ 
& ${\rm Ext}^1_{\mathcal{MM}}(H^1(A),M)=\frac{A(k)\otimes W_0M}
{{\rm Hom}_{\mathcal{MM}}(H^1(A),M/W_0M)}$\\
${\rm Ext}^{\ge 2}_{{\mathcal A}dm}(A(F)/A(k),W)=0$ & 
${\rm Ext}^{\ge 2}_{\mathcal{MM}}(H^1(A),M)=0$
\end{tabular} \end{proposition}
As $A(F)/A(k)$ is a canonical direct ``$H_1$''-summand of 
${\mathbb B}^{\bullet}(A)$, we see that admissible 
representations of finite length should be related to effective 
motives. At least the Ext's between some irreducible objects are dual. 

\subsection{${\mathcal I}_G$}
The formal properties of ${\mathcal A}dm$ are not very nice. In particular, 
to prove Theorem \ref{adm-ab} and Proposition \ref{ext-calc-ab} and to give 
an evidence to Conjecture \ref{conj-equi}, one uses the inclusion of 
${\mathcal A}dm$ to a bigger full subcategory in the category of smooth 
representations of $G$. 

\vspace{4mm}

{\sc Definition.} An object $W\in{\mathcal S}m_G$ is called ``homotopy 
invariant'' (in birational sense) if $W^{G_{F/L}}=W^{G_{F/L'}}$ for 
any purely transcendental subextension $L'/L$ in $F/k$. Denote by 
${\mathcal I}_G$ the full subcategory in ${\mathcal S}m_G$ with 
``homotopy invariant'' objects. 

\vspace{4mm}

{\sc Remark.} In this definition it suffices to 
consider only $L'$'s of finite type over $k$, cf. \cite{repr}, \S 6. 

\vspace{4mm}

A typical object of ${\mathcal I}_G$ is the ${\mathbb Q}$-space 
$CH^q(X_F)_{{\mathbb Q}}$ of cycles of codimension $q$ on 
$X\times_kF$ modulo rational equivalence, for any smooth variety 
$X$ over $k$. 
\begin{theorem}[\cite{repr}, $n=\infty$] \label{mnogo-o-I} 
\begin{enumerate} \item The category ${\mathcal I}_G$ 
is a Serre subcategory in ${\mathcal S}m_G$. 
\item \label{Adm-in-I} ${\mathcal A}dm\subset{\mathcal I}_G$, 
i.e. any admissible representation of $G$ is ``homotopy invariant''. 
\item \label{left-right-adj} The inclusion ${\mathcal I}_G
\hookrightarrow{\mathcal S}m_G$ admits a left and a right adjoints 
${\mathcal I},-^{(0)}:{\mathcal S}m_G\longrightarrow{\mathcal I}_G$. 
\item The objects $C_{k(X)}:={\mathcal I}{\mathbb Q}[\{k(X)
\stackrel{/k}{\hookrightarrow}F\}]$ for all irreducible varieties $X$ 
over $k$ form a system of projective generators of ${\mathcal I}_G$. 
\item \label{pse-mot-filt} For any smooth proper variety $X$ over $k$ 
there is a canonical filtration 
$C_{k(X)}\supset{\mathcal F}^1\supset{\mathcal F}^2\supset\dots$, 
canonical isomorphisms $C_{k(X)}/{\mathcal F}^1={\mathbb Q}$ and 
${\mathcal F}^1/{\mathcal F}^2={\rm Alb}(X_F)_{{\mathbb Q}}$, and a 
non-canonical splitting $C_{k(X)}\cong{\mathbb Q}\oplus
{\rm Alb}(X_F)_{{\mathbb Q}}\oplus{\mathcal F}^2$. The term 
${\mathcal F}^2$ is determined by these conditions together with 
${\rm Hom}_G({\mathcal F}^2,{\mathbb Q})={\rm Hom}_G
({\mathcal F}^2,A(F)/A(k))=0$ for any abelian variety $A$ over $k$. 
\item \label{qkobian} For any smooth proper variety $X$ over $k$ there 
is a canonical surjection $C_{k(X)}\to CH_0(X_F)_{{\mathbb Q}}$, which 
is injective if $X$ unirational over a curve.\footnote{and in some 
other cases when $CH_0(X)$ is ``known''} 
\item There exist (co-) limits in ${\mathcal I}_G$. \end{enumerate}
\end{theorem}

The following two conjectures link ${\mathcal I}_G$ 
with algebraic geometry and topology, 
\begin{conjecture}[\cite{repr}] \label{C=CH} If $n=\infty$ 
then the natural surjection $C_{k(X)}\longrightarrow 
CH_0(X\times_kF)_{{\mathbb Q}}$ is an isomorphism for 
any smooth proper variety $X$ over $k$. \end{conjecture} 

{\sc Remarks.} 1. One deduces from Theorem \ref{mnogo-o-I} 
(\ref{qkobian}) a description of the category of abelian varieties 
over $k$ with the groups of morphisms tensored with ${\mathbb Q}$ 
as a full subcategory of ${\mathcal A}dm_G\subset{\mathcal I}_G$ 
in terms of a functorial increasing ``level'' filtration 
$N_{\bullet}$ on smooth $G$-modules introduced in \cite{repr}. 

2. \label{ekviva-izo-vseh} The conjecture of Bloch and Beilinson 
(\cite{height} and \cite{bloch}, Lecture 1) on the ``motivic'' 
filtration on the Chow groups together with the semi-simplicity 
``standard'' conjecture of Grothendieck (asserting that numerical 
and homological equivalence coincide for smooth proper varieties), 
imply that ``numerical'' equivalence coincides with rational 
equivalence on the cycles on ${\bf Spec}$ of the tensor product of 
two fields over a common subfield (\cite{n-mot,mpi}). If combined 
with Conjecture \ref{C=CH}, this would give that 
${\mathbb B}^{\bullet}$ is an equivalence of categories (Conjecture 
\ref{conj-equi}), cf. also ``Corollary'' \ref{dva-sled} below.  

3. Define a binary non-associative operation 
$\otimes_{{\mathcal I}}$ on ${\mathcal S}m_G$ by $W_1
\otimes_{{\mathcal I}}W_2:={\mathcal I}(W_1\otimes W_2)$. 

It follows from Conjecture \ref{C=CH} that there is a canonical
isomorphism, the ``K\"unneth formula'': $C_{k(X\times_kY)}
\stackrel{\sim}{\longrightarrow}C_{k(X)}\otimes_{{\mathcal I}}
C_{k(Y)}$ for any pair of irreducible $k$-varieties $X,Y$. 
An evidence (and an inconditional proof in the case when $X$ 
is a curve) for this can be found in \cite{obz}. 

It would follow from the ``K\"unneth formula'' that the restriction 
of $\otimes_{{\mathcal I}}$ to ${\mathcal I}_G$ is a commutative 
associative tensor structure, and that the class of projective objects 
is stable under $\otimes_{{\mathcal I}}$, cf. \cite{repr}. 

It would be interesting to find a ``semi-simple graded'' version of 
$\otimes_{{\mathcal I}}$ to make ${\mathbb B}^{\bullet}$ a tensor functor. 

\begin{conjecture}[\cite{pgl}] \label{I-Omega} 
Any irreducible object of ${\mathcal I}_G$ is contained 
in the algebra $\Omega^{\bullet}_{F/k}$ if $n=\infty$. \end{conjecture} 

\begin{conj-corollary}[\cite{pgl}] \label{dva-sled} \begin{itemize} 
\item If numerical equivalence coincides with homological 
then ${\mathbb B}^{\bullet}$ is an equivalence of categories. 
\item Any irreducible object of ${\mathcal I}_G$ 
is admissible if $n=\infty$. So ``${\mathcal I}_G\approx{\mathcal A}dm$''. 
\end{itemize} \end{conj-corollary}

Conjecture \ref{I-Omega} is one of the main motivations for the study 
of semi-linear representations of $G$, cf. \S\ref{ot-lin-k-polulin}. 
It has also the following geometric corollary, conjectured by Bloch. 
\begin{conj-corollary}[\cite{pgl}] \label{gip-bloha} If 
$\Gamma(X,\Omega^{\ge 2}_{X/k})=0$ for a smooth proper variety 
$X$ over $k$ then the Albanese map induces an isomorphism 
$CH_0(X)^0\stackrel{\sim}{\longrightarrow}{\rm Alb}(X)$. 
In that case $C_{k(X)}=CH_0(X_F)_{{\mathbb Q}}$. \end{conj-corollary}

{\sc Remarks.} 1. There is a locally compact group $H$ and 
a continuous injective homomorphism with dense image 
$H\longrightarrow G$ such that ${\mathcal I}_G$ admits an 
explicit description as a full subcategory of ${\mathcal S}m_H$ 
stable under taking subquotients (but not extensions). 

The category ${\mathcal S}m_H$ may be useful in the 
study of left derivatives of additive functors. 
 
2. ${\mathcal I}_G$ is equivalent to the category of non-degenerate 
modules over an associative idempotented algebra, \cite{obz}. 

\subsection{Differential forms}
In an attempt to compare various cohomology theories $H^{\ast}$, one 
can associate with them some $G$-modules, like $H^{\ast}(F):=\lim\limits
_{\longrightarrow}H^{\ast}(U)$, where $U$ runs over spectra of smooth 
subalgebras in $F$ of finite type over $k$, or the image $H^{\ast}_c(F)$ 
in $H^{\ast}(F)$ of $\lim\limits_{\longrightarrow}H^{\ast}(X)$, where $X$ 
runs over smooth proper models of subfields in $F$ of finite type over $k$. 

Clearly, $H^{\ast}_c(F)$ is an admissible representation of $G$ over 
$H^{\ast}(k)$. It would follow from the semi-simplicity standard 
conjecture that it is semi-simple. For instance, it replaces reference 
to the semi-simplicity 
standard conjecture in Remark 2 on p.\pageref{ekviva-izo-vseh}. 

In the case $H^{\ast}=H^{\ast}_{{\rm dR}/k}$ of the de Rham cohomology 
the graded quotients of the (descending) Hodge filtration on 
$H^q_{{\rm dR}/k,c}(F)$ are $H^{p,q-p}_{F/k}=\lim\limits_{\longrightarrow}
{\rm coker}[H^{p-1}(D,\Omega^{q-p-1}_{D/k})\longrightarrow 
H^p(X,\Omega^{q-p}_{X/k})]$, where $(X,D)$ runs over 
pairs consisting of a smooth proper variety $X$ with $k(X)\subset F$ 
and a normal crossing divisor $D$ on $X$ with smooth irreducible 
components. More particularly, $H^{q,0}_{F/k}=\Omega^q_{F/k,{\rm reg}}
\subset H^q_{{\rm dR}/k,c}(F)$ is the $G$-submodule spanned by the 
spaces $\Gamma(X,\Omega^{\bullet}_{X/k})$ of regular differential 
forms on all smooth projective $k$-varieties $X$ with the function 
fields embedded into $F$. 

\begin{proposition}[\cite{obz}] \label{reg-top-semi} Suppose that the 
cardinality of $k$ is at most continuum. Fix an embedding $\iota:k
\hookrightarrow{\mathbb C}$ to the field of complex numbers. Then 
\begin{itemize} \item there is a ${\mathbb C}$-anti-linear canonical 
isomorphism (depending on $\iota$) $H^{p,q}_{F/k}\otimes_{k,\iota}
{\mathbb C}\cong H^{q,p}_{F/k}\otimes_{k,\iota}{\mathbb C}$; 
\item the representation $H^n_{{\rm dR}/k,c}(F)$ (and thus, 
$\Omega^n_{F/k,{\rm reg}}$) is semi-simple for any $1\le n<\infty$. 
\end{itemize} \end{proposition}

Recall (Theorem \ref{mnogo-o-I}(\ref{left-right-adj})), that for any 
$W\in{\mathcal S}m_G$ its maximal subobject of $W$ in ${\mathcal I}_G$, 
the ``homotopy invariant'' part of $W$, is denoted by $W^{(0)}$. The 
following fact gives one more evidence for the cohomological nature of 
the objects of ${\mathcal I}_G$, since $\Omega^{\bullet}_{F/k,{\rm reg}}$ 
is the ``cohomological part'' of $\bigotimes^{\bullet}_F\Omega^1_{F/k}$. 
\begin{proposition}[\cite{pgl}] 
If $n=\infty$ then $(\bigotimes^{\bullet}_F\Omega^1_{F/k})^{(0)}
=\Omega^{\bullet}_{F/k,{\rm reg}}$. \end{proposition} 

\begin{proposition}[\cite{obz}] \label{1-forms} For any 
$1\le n\le\infty$ the representation $H^1_{{\rm dR}/k}(F)$ modulo 
the sum of submodules isomorphic to $F^{\times}/k^{\times}$ is a 
direct sum of $\# k$ copies of $A(F)/A(k)$ for all isogeny classes 
$A$ of simple abelian $k$-varieties. In particular, 
$\Omega^1_{F/k,{\rm reg}}$ is semi-simple. \end{proposition} 

This suggests that the isomorphism classes of irreducible 
subquotients of $H^{\ast}_c(F)$ can be naturally identified with the 
irreducible effective primitive motives, and that the isomorphism 
classes of irreducible subquotients of $H^{\ast}(F)$ are 
related to more general irreducible effective motives, such as the 
Tate motive ${\mathbb Q}(-1)$ in the case of $H^1_{{\rm dR}/k}(F)$. 

\section{From linear to semi-linear representations} 
\label{ot-lin-k-polulin}
The representation $\Omega^{\bullet}_{F/k}$ of $G$ is also 
an $F$-vector space endowed with a semi-linear $G$-action. 

\vspace{4mm}

{\sc Definition.} A semi-linear representation of $G$ over $F$ is an 
$F$-vector space $V$ endowed with an additive $G$-action $G\times V\to V$ 
such that $g(fv)=gf\cdot gv$ for any $g\in G$, $v\in V$ and $f\in F$. 

\vspace{4mm}

Denote by ${\mathcal C}$ the category of smooth semi-linear 
representations of $G$ over $F$. 

It is well-known after Hilbert, Tate, Sen, Fontaine... that the 
semi-linear representations is a powerful tool in the study of 
Galois representations. We try to use them in non-Galois context.  

In some respects ${\mathcal C}$ is simpler than ${\mathcal S}m_G$. 
In particular, it follows from Hilbert's Theorem 90 that the category 
${\mathcal C}$ admits a countable system of cyclic generators: 
$F[G/G_{F/K_m}]$, where $K_m$ is a purely transcendental extension 
of $k$ in $F$ of transcendence degree $m$. 

\vspace{4mm}

Once again, we are interested in linear representations of $G$, 
especially in irreducible ones, and more particularly, in irreducible 
``homotopy invariant'' representations, i.e. objects of ${\mathcal I}_G$. 

The problem of describing (the irreducible objects of) ${\mathcal S}m_G$ 
could be split into describing (the irreducible objects of) ${\mathcal C}$ 
and their linear submodules. 

For example, all representations $A(F)/A(k)$ 
of $G$ for all abelian $k$-varieties $A$ (i.e. corresponding to all pure 
1-motives) are contained in the irreducible object $\Omega^1_{F/k}$ of 
${\mathcal C}$. 

\vspace{4mm}

Suppose from now on that $n=\infty$. 

\vspace{4mm}

One has the faithful forgetful functor ${\mathcal C}
\stackrel{{\rm for}}{\longrightarrow}{\mathcal S}m_G(k)$ 
admitting a left adjoint functor of extending of coefficients to $F$: 
${\mathcal S}m_G(k)\stackrel{\otimes_kF}{\longrightarrow}{\mathcal C}$, 
where ${\mathcal S}m_G(k)$ is the category of smooth representations 
of $G$ over $k$, so $W\hookrightarrow{\rm for}(W\otimes_kF)$. The 
functor $F\otimes_k$ is not full and does not respect the irreducibility. 

However, if $W$ is irreducible, there is an irreducible semi-linear 
quotient $V$ of $W\otimes F$ with an inclusion $W\subset V$, so any 
irreducible object of ${\mathcal S}m_G$ is contained in an 
irreducible object of ${\mathcal C}$. 

This gives a hint that it 
might be sufficient for the study of some categories of 
${\mathbb Q}$-linear representations of $G$ to know the structure of 
some ``relatively small'' full sub-category of ${\mathcal C}$. 

\vspace{4mm}

The following claim suggests the category ${\mathcal C}$ is 
``more complicated'' than ${\mathcal I}_G(k)$. However, this should 
be compaired with Lemma \ref{Glob-kvazi-obr}. 
\begin{lemma}[\cite{adm}, Lemma 0.1] \label{exe-odin-strog} 
The functor ${\mathcal I}_G(k)\stackrel{F\otimes_k}
{\longrightarrow}{\mathcal C}$ is fully faithful. \end{lemma}

Another, though a weaker, but a little bit more explicit 
condition on the semi-linear quotients of $W\otimes F$ for 
$W\in{\mathcal I}_G$ is given in the next section \ref{usl-glob-porozhd}. 

\subsection{Valuations and associated functors (\cite{max})}
\label{usl-glob-porozhd}
In order to associate functors on categories of $k$-varieties to 
representations of $G$ one can try to ``approximate'' rings by their 
subfields. Evidently, this does not work literally, but apparently 
works in the case of discrete valuation rings of $F$. 

Let $v:F^{\times}/k^{\times}\longrightarrow{\mathbb Q}$ be a discrete 
valuation, ${\mathcal O}_v$ be the valuation ring, $\mathfrak{m}_v
={\mathcal O}_v-{\mathcal O}_v^{\times}$ be the maximal ideal, and 
$\kappa(v)$ be the residue field. Denote by 
${\mathcal P}_F$ the set of all such valuations. 

Set $G_v:=\{\sigma\in G~|~\sigma({\mathcal O}_v)=
{\mathcal O}_v\}$. This is a closed subgroup in $G$. 
The $G_v$-action on $\kappa(v)$ induces a homomorphism 
$G_v\longrightarrow\hspace{-3mm}\to G_{\kappa(v)/k}$. 

\begin{proposition} For any discrete valuation 
$v\in{\mathcal P}_F$ the additive functor $(-)_v:{\mathcal S}m_G
\longrightarrow{\mathcal S}m_{G_v}$, $W\mapsto W_v:=
\sum_{F'\subset{\mathcal O}_v}W^{G_{F/F'}}\subseteq W$, is fully 
faithful and preserves surjections and injections. \end{proposition}

Then the additive subfunctor $\Gamma:{\mathcal S}m_G\longrightarrow
{\mathcal S}m_G$ of the identity functor, defined by $W\mapsto
\Gamma(W):=\bigcap_{v\in{\mathcal P}_F}W_v$, preserves injections. 

{\sc Example.} \label{Gamma-ot-1-form} $\Gamma(\Omega^1_{F/k})
\cong\bigoplus_A(A(F)/A(k))\otimes_{{\rm End}A}\Gamma(A,\Omega^1_{A/k})$, 
where $A$ runs over the set of isogeny classes of simple abelian 
varieties over $k$. 
\begin{lemma} \label{Glob-kvazi-obr} The compositions 
${\mathcal I}_G(k)\stackrel{F\otimes_k}{\longrightarrow}{\mathcal C}
\stackrel{\Gamma\circ{\rm for}}{\longrightarrow}{\mathcal S}m_G(k)$ 
and ${\mathcal I}_G\hookrightarrow{\mathcal S}m_G\stackrel{\Gamma}
{\longrightarrow}{\mathcal S}m_G$ are identical. \end{lemma}

{\sc Remark.} This implies that any semi-linear quotient $V$ of 
$W\otimes F$ with $W\in{\mathcal I}_G$ (in particular, any irreducible 
semi-linear representation $V$ containing a ``homotopy invariant'' 
representation), is ``globally generated'', i.e., $\Gamma(V)\otimes F
\longrightarrow\hspace{-3mm}\to V$ is surjective. 

This is the condition one can impose on the class of ``interesting'' 
semi-linear representations. There are some reasons to expect that 
$(-)_v$ is exact, cf. \cite{max}. This would imply some nice 
properties of the category of ``globally generated'' semi-linear 
representations. 

\subsection{Admissible semi-linear representations}
In the study of representations of any group, it is natural 
to start with the finite-dimensional representations. 

\begin{theorem}[\cite{repr}] Any finite-dimensional smooth semi-linear 
representation of $G$ over $F$ is trivial, if $n=\infty$.\end{theorem}

\vspace{4mm}

A natural extension of the notion of finite-dimensional semi-linear 
representation in the case $n=\infty$ is the notion of admissible 
semi-linear representation. 

{\sc Definition.} A smooth semi-linear representation $V$ of $G$ over $F$ 
is called {\sl admissible} if, for any open subgroup $U\subseteq G$, the 
fixed subspace $V^U$ is finite-dimensional over the fixed subfield $F^U$ 
(or equivalently, $\dim_LV^{G_{F/L}}<\infty$ for any subfield 
$L\subset F$ of finite type over $k$). 

\begin{theorem}[\cite{pgl,adm}] The admissible semi-linear 
representations of $G$ over $F$ form an abelian tensor 
(but not rigid) category, denoted by ${\mathcal A}$. 

The functor 
$H^0(G_{F/L},-)$ is exact on ${\mathcal A}$ for any subfield 
$L\subseteq F$, so $F$ is a projective object of ${\mathcal A}$. 
\end{theorem} 

\vspace{4mm}
The latter give an example of an admissible semi-linear representation. 

{\sc Example.} Let the ideal ${\mathfrak m}\subset F
\otimes_{\overline{{\mathbb Q}}}F$ be the kernel of the multiplication 
map $F\otimes_{\overline{{\mathbb Q}}}F\stackrel{\times}{\longrightarrow}F$. 
Consider the powers of the ideal ${\mathfrak m}$ as objects of ${\mathcal C}$ 
with the $F$-multiplication via $F\otimes_{\overline{{\mathbb Q}}}
\overline{{\mathbb Q}}$. Then ${\mathfrak m}^s/{\mathfrak m}^{s+1}
={\rm Sym}^s_F\Omega^1_F$, and the objects ${\rm Sym}^s_F
\Omega^{\bullet}_{F/k}$ are admissible for all $s\ge 1$. 

\vspace{4mm}

In the case $k=\overline{{\mathbb Q}}$, the field of algebraic 
complex numbers, the category ${\mathcal A}$ is equivalent to the 
category of ``coherent'' sheaves in smooth topology $\mathfrak{Sm}_k$ 
on $k$. (The underlying category of $\mathfrak{Sm}_k$ is the smooth 
morphisms of smooth $k$-varieties, and the coverings are coverings 
by images; ``coherent'' means: a sheaf of ${\mathcal O}$-modules such 
that its restriction to the small Zariski site of any smooth $k$-variety 
is coherent.) Moreover, ${\mathcal A}$ admits the following explicit 
description, cf. \cite{adm}. 
\begin{itemize} \item The sum of the images of the $F$-tensor powers 
$\bigotimes^{\ge\bullet}_F{\mathfrak m}$ under all morphisms in 
${\mathcal C}$ defines a decreasing filtration $W^{\bullet}$ on the 
objects of ${\mathcal A}$ such that its graded quotients $gr^q_W$ are 
finite direct sums of direct summands of $\bigotimes^q_F\Omega^1_F$. 
This filtration is evidently functorial and multiplicative: 
$(W^pV_1)\otimes_F(W^qV_2)\subseteq W^{p+q}(V_1\otimes_FV_2)$ 
for any $p,q\ge 0$ and any $V_1,V_2\in{\mathcal A}$. 
\item ${\mathcal A}$ is equivalent to the direct sum of the category 
of finite-dimensional $k$-vector spaces and its abelian full 
subcategory ${\mathcal A}^{\circ}$ with objects $V$ such that $V^G=0$. 
\item Any object $V$ of ${\mathcal A}^{\circ}$ is a quotient of 
a direct sum of objects (of finite length) of type $\bigotimes^q_F
({\mathfrak m}/{\mathfrak m}^s)$ for some $q,s\ge 1$. 
\item If $V\in{\mathcal A}$ is of finite type then it is 
of finite length and $\dim_k{\rm Ext}_{{\mathcal A}}^j(V,V')
<\infty$ for any $j\ge 0$ and any $V'\in{\mathcal A}$; 
if $V\in{\mathcal A}$ is irreducible and ${\rm Ext}^1_{{\mathcal A}}
({\mathfrak m}/{\mathfrak m}^q,V)\neq 0$ for some $q\ge 2$ then 
$V\cong{\rm Sym}^q_F\Omega^1_F$ and ${\rm Ext}^1_{{\mathcal A}}
({\mathfrak m}/{\mathfrak m}^q,V)\cong k$. 
\item ${\mathcal A}^{\circ}$ has no projective objects, but 
$\bigotimes^q_F{\mathfrak m}$ are its ``projective pro-generators'': 
the functor ${\rm Hom}_{{\mathcal C}}(\bigotimes^q_F{\mathfrak m},-)
=\lim\limits_{\longrightarrow}{\rm Hom}_{{\mathcal A}}
(\bigotimes^q_F({\mathfrak m}/{\mathfrak m}^N),-)$ 
is exact on ${\mathcal A}$ for any $q$. 
\end{itemize}

\vspace{4mm}

Representations of particular interest are admissible ones. 
Though tensoring with $F$ does not transform them to 
admissible semi-linear representations, there exists a similar 
functor in the opposite direction, faithful 
at least if $k=\overline{{\mathbb Q}}$,.

It is explained in \cite{adm}, that when $k=\overline{{\mathbb Q}}$, 
for any object $V$ of ${\mathcal A}$ and any smooth $k$-variety $Y$, 
embedding of the generic points of $Y$ into $F$ determines a locally 
free coherent sheaf ${\mathcal V}_Y$ on $Y$ with the generic fibre 
$V^{G_{F/k(Y)}}$. Moreover, for any dominant morphism $X\stackrel
{\pi}{\longrightarrow}Y$ of smooth $k$-varieties, the inclusion
of the generic fibres $k(X)\otimes_{k(Y)}V^{G_{F/k(Y)}}\subseteq
V^{G_{F/k(X)}}$ induces an injection of the coherent sheaves 
$\pi^{\ast}{\mathcal V}_Y\hookrightarrow{\mathcal V}_X$ on $X$, 
which is an isomorphism if $\pi$ is {\'e}tale.

For any $V\in{\mathcal A}$ the space $\Gamma(Y,{\mathcal V}_Y)$ 
is a birational invariant of smooth proper $Y$. Then we get a 
left exact functor ${\mathcal A}\stackrel{\Gamma}{\longrightarrow}
{\mathcal S}m_G(k)$ given by $V\mapsto\lim\limits_{\longrightarrow}
\Gamma(Y,{\mathcal V}_Y)$, where $Y$ runs over the smooth proper 
models of subfields in $F$ of finite type over $k$. In general, 
$\Gamma(V)$ is not admissible.

The functor $\Gamma$ coincides with the composition of the forgetful 
functor to the category of smooth representations of $G$ with the 
functor $\Gamma$ from \S\ref{usl-glob-porozhd}. The functor $\Gamma$ 
is faithful, but it is not full, and the objects in its image are 
highly reducible, cf. Example on p.\pageref{Gamma-ot-1-form}. 

\begin{conjecture} \label{ustr-dop} \begin{enumerate} 
\item The functor ${\rm Hom}_{{\mathcal C}}(\otimes^q_F{\mathfrak m},-)$ 
is exact on ${\mathcal A}$ for any $q\ge 0$. 
\item \label{part-2} Irreducible objects of ${\mathcal A}$ are direct 
summands of the tensor algebra $\bigotimes^{\bullet}_F\Omega^1_{F/k}$. 
\item ${\mathcal A}$ is equivalent to the category of
``coherent'' sheaves on $\mathfrak{Sm}_k$.
\end{enumerate}\end{conjecture}
As another evidence for Conjecture \ref{ustr-dop} (\ref{part-2}), 
in addition to the case $k=\overline{{\mathbb Q}}$, it is shown in 
\cite{pgl} that for any $L\subset F$ purely transcendental of degree 
$m$ over $k$ and any $V\in{\mathcal A}$ any irreducible subquotient of 
the $L$-semi-linear representation $V^{G_{F/L}}$ of ${\rm PGL}_{m+1}k$ 
is a direct summand of $\bigotimes_L^{\bullet}\Omega^1_{L/k}$. 

As there exist smooth non-admissible irreducible semi-linear 
representations, cf. \cite{pgl}, \S4.2, one cannot replace the category 
${\mathcal A}$ in the part (\ref{part-2}) of Conjecture \ref{ustr-dop} 
by the whole category ${\mathcal C}$, and has to put some additional 
conditions, e.g. the one mentioned in \S\ref{usl-glob-porozhd}. 

{\sc Remark.} Assuming the part (\ref{part-2}) of Conjecture \ref{ustr-dop}, 
one can reformulate Conjecture \ref{I-Omega} in the following 
linguistically more convencing form: 

{\it Any irreducible object of ${\mathcal A}dm$ (and of ${\mathcal I}_G$) 
is contained in an irreducible object of ${\mathcal A}$.}

\vspace{5mm}

\noindent
{\sl Acknowledgement.} {\small I would like to thank Uwe Jannsen 
for his encouraging interest in this work, many suggestions and 
many inspiring discussions. 

This paper is an expanded write-up of several seminar talks given 
at Steklov Institute in Moscow, at G\"{o}ttingen University and 
at the Max-Plank-Institut in Bonn. I would like to thank Dmitry Orlov, 
Yuri Tschinkel and Yu.I.Manin for inviting me to these seminars. 

I am grateful to the Max-Planck-Institut f\"ur Mathematik in Bonn 
and to Regensburg University for their hospitality. 

In the course of this work the author was supported by the following 
institutions: the Max-Planck-Institut f\"ur Mathematik in Bonn 
(September--November 2004), Alexander von Humboldt-Stiftung 
(March--July 2005 and September 2005--March 2006) and Russian 
Foundation for Basic Research (under grant 02-01-22005). 
 
}

\newpage
\appendix
\section{A construction of a separable closure of a transcendence degree 
one extension of an algebraically closed field of positive characteristic}
The following claim is a refinement of Proposition 4.1 of \cite{repr}. 
\begin{proposition} \label{add-mult} The $G^{\circ}$-modules 
$F/k$ and $F^{\times}/k^{\times}$ are irreducible if either 
${\rm char}(k)=0$, or $2\le n\le\infty$. If $n=1$ and ${\rm char}
(k)\neq 0$ then the $G^{\circ}$-orbit of $x$ generates the separable 
closure of $k(x)$ in $F$ for any $x\in F$. \end{proposition} 
{\it Proof.} Let $A$ be the additive subgroup of $F$ generated by 
the $G^{\circ}$-orbit of some $x\in F-k$. For any $y\in A-k$ one has 
$\frac{2}{y^2-1}=\frac{1}{y-1}-\frac{1}{y+1}$. As $\frac{1}{y-1}$ 
and $\frac{1}{y+1}$ are in the $G^{\circ}$-orbit of $y$, this 
implies that $y^2\in A$. As for any $y,z\in A$ one has 
$yz=\frac{1}{4}((y+z)^2-(y-z)^2)$, the group $A$ is 
a subring of $F$, if ${\rm char}(k)\neq 2$. 

Let $M$ be the multiplicative subgroup of $F^{\times}$ generated by 
the $G^{\circ}$-orbit of some $x\in F-k$. Then for any $y,z\in M$ 
one has $y+z=z(y/z+1)$, so if $y/z\not\in k$ then $y+z\in M$, and 
thus, $M\bigcup\{0\}$ is a $G^{\circ}$-invariant subring of $F$. 

Since the $G^{\circ}$-orbit of an element $x\in F-k$ contains all 
elements of $F-\overline{k(x)}$, if $n\ge 2$ then each element 
of $F$ is the sum of a pair of elements in the orbit. Any 
$G^{\circ}$-invariant subring in $F$, but not in $k$, is a 
$k$-subalgebra, so if $n=1$ then ${\rm Gal}(F/{\mathbb Q}(G^{\circ}x))
\subset G^{\circ}$ is a compact subgroup normalized by 
$G^{\circ}$. Then by Theorem 2.9 of \cite{repr} we have 
${\rm Gal}(F/{\mathbb Q}(G^{\circ}x))=\{1\}$, i.e., the extension 
$F/{\mathbb Q}(G^{\circ}x)$ is purely inseparable. As any element 
of ${\mathbb Q}(G^{\circ}x)$ is the fraction of a pair of 
elements in ${\mathbb Z}[G^{\circ}x]$ and for any $y\in F-k$ 
the element $1/y$ belongs to the $G^{\circ}$-orbit of $y$, the 
${\mathbb Z}$-subalgebra generated by the $G^{\circ}$-orbit of 
$x$ coincides with $F$, if ${\rm char}(k)=0$, or $2\le n\le\infty$. 

Let us show that $k(G^{\circ}x)$ is a separable extension of 
$k(x)$, equivalently, that if $\sigma^Nx=x$ for some $N\ge 1$ 
then $k(x,\sigma x)$ is a separable extension of $k(x)$. Let 
$P(x,\sigma x)$ be a minimal polynomial. Then $P_Idx+P_{II}
d(\sigma x)=0\in\Omega^1_{k(x,\sigma x)/k}$, where either 
$P_I\neq 0$, or $P_{II}\neq 0$ as otherwise $P=Q^p$ for 
another polynomial $Q$. If $P_{II}\neq 0$ then $k(x,\sigma x)$ 
is a separable extension of $k(x)$. If $P_I\neq 0$ then 
$k(x,\sigma x)$ is a separable extension of $k(\sigma x)$, and 
thus, $k(x,\sigma^{-1}x)$ is a separable extension of $k(x)$. 
Then $k(x,\sigma^{-1}x,\dots,\sigma^{-(N-1)}x=\sigma x)$ 
is a separable extension of $k(x)$. \qed 

\section{The ``K\"{u}nneth formula'' for products with curves}
Define a $G$-homomorphism $${\mathbb Q}[\{k(X)\otimes_kk(Y)
\stackrel{/k}{\hookrightarrow}F\}]\stackrel{\alpha}{\longrightarrow}
C_{k(X)}\otimes C_{k(Y)}\quad\mbox{by}\quad\tau\mapsto\tau|_{k(X)}
\otimes\tau|_{k(Y)}.$$ It is shown in \cite{repr} that $\alpha$ 
is surjective, which gives a surjection $C_{k(X\times_kY)}
\longrightarrow C_{k(X)}\otimes_{{\mathcal I}}C_{k(Y)}$. 

For arbitrary $A\in C_{k(X)}$ and $B\in C_{k(Y)}$ choose some liftings 
$\widetilde{A}\in{\mathbb Q}[\{k(X)\stackrel{/k}
{\hookrightarrow}F\}]$ and $\widetilde{B}\in{\mathbb Q}
[\{k(Y)\stackrel{/k}{\hookrightarrow}F\}]$ such that 
all embeddings from $\widetilde{A}$ and from $\widetilde{B}$ 
are is pairwise general position.\footnote{First, choose arbitrary 
$\widetilde{A}$ and $\widetilde{B}$. For each point $P$ of the support 
of $\widetilde{B}$ choose a generic curve $C$ passing through $P$, on 
which $P$ is a generic point with respect to a field of definition of 
$C$. Replace $P$ by a linearly equivalent linear combination of points 
of $C$ in general position with respect to $\widetilde{A}$. Then 
we get the desired $\widetilde{B}$.}

One has to check that the class of $\widetilde{A}\times\widetilde{B}
\in{\mathbb Q}[\{k(X\times_kY)\stackrel{/k}{\hookrightarrow}F\}]$ 
in $C_{k(X\times_kY)}$ is independent of the choice of $\widetilde{A}$ 
and $\widetilde{B}$. If some other liftings 
$\widetilde{A}'\in{\mathbb Q}[\{k(X)\stackrel{/k}{\hookrightarrow}
F\}]$ and $\widetilde{B}'\in{\mathbb Q}[\{k(Y)\stackrel{/k}
{\hookrightarrow}F\}]$ are defined similarly, choose some lifting 
$\widetilde{B}''\in{\mathbb Q}[\{k(Y)\stackrel{/k}{\hookrightarrow}F\}]$ 
of $B$ such that all embeddings from $\widetilde{A}$ and from 
$\widetilde{B}''$, as well as from $\widetilde{A}'$ and from 
$\widetilde{B}''$, are in pairwise general position. Then 
$\widetilde{A}\times\widetilde{B}-\widetilde{A}'\times
\widetilde{B}'=(\widetilde{A}-\widetilde{A}')\times\widetilde{B}''
+\widetilde{A}\times(\widetilde{B}-\widetilde{B}'')+
\widetilde{A}'\times(\widetilde{B}''-\widetilde{B}')$.

Thus, one has to check the following condition $\star_{X,Y}$: 
if the class of $\sum\limits_{i=1}^Na_i\tau_i\in{\mathbb Q}[\{k(X)
\stackrel{/k}{\hookrightarrow}F\}]$ in $C_{k(X)}$ is zero and all 
$\tau_i$ are in general position with respect to $\sigma:
k(Y)\stackrel{/k}{\hookrightarrow}F$ then the class of $\gamma:=
\sum\limits_{i=1}^Na_i(\tau_i,\sigma)\in{\mathbb Q}[\{k(X\times_kY)
\stackrel{/k}{\hookrightarrow}F\}]$ in $C_{k(X\times_kY)}$ is zero. 
Also, one has to check the condition $\star_{Y,X}$. 

By definition of the functor ${\mathcal I}$, there exist purely 
transcendental extensions $L_j'/L_j$, elements $\alpha_j\in
{\mathbb Q}[\{k(X)\stackrel{/k}{\hookrightarrow}F\}]^{G_{F/L_j'}}$
and $\xi_j\in G_{F/L_j}$ such that $\sum\limits_{i=1}^Na_i\tau_i
=\sum_j(\xi_j\alpha_j-\alpha_j)$.

If $\sigma$ is in general position with respect to the compositum 
$L$ of all $\tau_i(k(X))$ then there exists $\kappa\in G_{F/L}$ such 
that $\kappa\sigma=:\sigma'$ is in general position with respect to 
the compositum of all $L_j'$. Then $\gamma':=\kappa\gamma=\sum_ia_i
(\tau_i,\sigma')=\sum_j(\xi_j\alpha_j-\alpha_j)\otimes\sigma'$.
Set $K_j:=L_j\sigma'(k(Y))$ and $K_j':=L_j'\sigma'(k(Y))$. Then 
$\alpha_j\otimes\sigma'\in{\mathbb Q}[\{k(X\times_kY)\stackrel{/k}
{\hookrightarrow}F\}]^{G_{F/K_j'}}$, $K_j'$ is a purely transcendental 
extension of  $K_j$, and there exist $\xi_j'\in G_{F/\sigma'(k(Y))}$ 
such that $\xi_j'|_{L_j'}=\xi_j|_{L_j'}$. This implies that $\gamma'=
\sum_j(\xi'_j(\alpha_j\otimes\sigma')-\alpha_j\otimes\sigma')$
belongs, by definition of the functor ${\mathcal I}$, to the kernel 
of the projection ${\mathbb Q}[\{k(X\times_kY)\stackrel{/k}
{\hookrightarrow}F\}]\longrightarrow C_{k(X\times_kY)}$, 
and therefore, the same is true for $\gamma$.

Let us check that the conditions $\star_{X,Y}$ and $\star_{Y,X}$
are equivalent. Consider a generic curve $C$ on $Y$, passing 
through $\sigma$, defined over a field containing the compositum 
of all $\tau_i(k(X))$. Then $\sigma$ is linearly equivalent to a 
linear combination $\beta$ of generic points of $C$ (which are 
therefore generic points of $Y$). Then the image of $\gamma$ 
in $C_{k(X\times_kY)}$ coincides with the image of 
$\sum_ia_i\tau_i\times(\sigma-\beta)$, which shows 
the implication $\star_{Y,X}\Rightarrow\star_{X,Y}$.

{\sc Example.} Let us check the condition $\star_{X,Y}$ in the case, 
when $X$ is a smooth proper curve. Let $K=\overline{\sigma(k(Y))}$. 
Then $\sum_ia_i\tau_i$ is a generic divisor on the curve $X_K$ over $K$, 
linearly equivalent to zero. According to Lemma 6.18 from \cite{repr}, 
the $G_{F/K}$-module of generic divisors on $X_K$ over $K$, linearly 
equivalent to zero, is generated by the elements $w_M=\sum\limits_{j=1}^M
(\sigma_j-\sigma_j')$ for all $M\gg 0$, where $(\sigma_1,\dots,
\sigma_M;\sigma_1',\dots,\sigma_M')$ is a generic $F$-point of the fibre 
over 0 of the morphism $X^M_K\times_KX^M_K\longrightarrow{\rm Pic}^0X_K$, 
sending $(x_1,\dots,x_M;y_1,\dots,y_M)$ to the class of 
$\sum\limits_{j=1}^M(x_j-y_j)$. Clearly, the compositum of all 
$\sigma_j(k(X))\sigma_j'(k(X))$ is in general position with respect to 
$K$. The same is true for any other element in the $G_{F/K}$-orbit of 
$w_M$. Therefore, as we have already seen above, the image of 
$\sum_ia_i(\tau_i,\sigma)$ in $C_{k(X\times_kY)}$ is zero. 

\vspace{4mm}

Thus, one has a canonical $G$-module surjection $C_{k(X)}\otimes 
C_{k(Y)}\longrightarrow\hspace{-3mm}\rightarrow C_{k(X\times_kY)}$, 
at least if $X$ is a curve,\footnote{and also if 
$C_{k(X)}=CH_0(X_F)_{{\mathbb Q}}$ and the transcendence degree of 
$k$ is infinite (or if the same is true for algebraic closures 
of all extensions of $k$ of finite type): $\sum_ia_i\tau_i
\in{\mathbb Q}[\{K(X)\stackrel{/K}{\hookrightarrow}F\}]$ is rationally 
equivalent to zero; as ${\rm tr.deg}(k)=\infty$, one has 
${\mathcal I}_{/K}{\mathbb Q}[\{K(X)\stackrel{/K}{\hookrightarrow}
F\}]=CH_0(X_F)_{{\mathbb Q}}$ (identify $k$ with $K$ by an automorphism 
of $F$ identical on the field of definition of $X$), i.e. 
$\sum_ia_i(\tau_i,\sigma)$ belongs to the kernel of the composition 
${\mathbb Q}[\{K(X)\stackrel{/K}{\hookrightarrow}F\}]
\stackrel{\sigma}{\hookrightarrow}{\mathbb Q}[\{k(X\times_kY)
\stackrel{/k}{\hookrightarrow}F\}]\longrightarrow\hspace{-2mm}
\rightarrow{\mathcal I}_{/K}{\mathbb Q}[\{k(X\times_kY)
\stackrel{/k}{\hookrightarrow}F\}]\longrightarrow\hspace{-2mm}
\rightarrow C_{k(X\times_kY)}.$} and the composition 
$C_{k(X\times_kY)}\longrightarrow C_{k(X)}\otimes_{{\mathcal I}}
C_{k(Y)}\longrightarrow C_{k(X\times_kY)}$ is identical. 

\section{Differential forms}
\begin{proposition} Let $k_0\subseteq k$ be a subfield. 
If $i<n$ then ${\rm Hom}_G(\Omega^i_F,\Omega^j_{F/k_0})
=\Omega^{j-i}_{k/k_0}\oplus\Omega^{j-i-1}_{k/k_0}\cdot d$; 
${\rm Hom}_G(\Omega^n_F,F)=0$ for any $n\ge 1$; ${\rm Hom}_G
(\Omega^q_{F/k_0,\log},\Omega^{\bullet}_{F/k_0})=\Omega^{\bullet-q}
_{k/k_0}$; ${\rm End}_G(\Omega^q_{F/k_0,\log})={\mathbb Q}$ 
for any $1\le q\le n$. \end{proposition}\vspace{-3mm}
{\it Proof.} For $i<n$ the $G$-module $\Omega^i_F$ is cyclic 
and generated by $\eta:=x_0dx_1\wedge\dots\wedge dx_i$ for some 
algebraically independent elements $x_0,\dots,x_i\in F-\overline{k}$. 
The space ${\mathbb Q}\cdot\eta\subset\Omega^i_F$ is stable 
under the subgroup $H:=G_{\{F,{\mathbb Q}^{\times}x_0,
{\mathbb Q}^{\times}x_1+{\mathbb Q},\dots,{\mathbb Q}^{\times}x_i+
{\mathbb Q}\}/k}\subset G$ and $H$ acts via its quotient 
$({\mathbb Q}^{\times})^{i+1}$. For any $\varphi\in{\rm Hom}_G
(\Omega^i_F,\Omega^j_{F/k_0})$ the form $\omega=\varphi(\eta)$ 
is $G_{\{F,x_1+{\mathbb Q},\dots,x_i+{\mathbb Q}\}/k(x_0)}$-invariant, 
so $\omega=\sum\limits_{S=(0\le s_1<\dots<s_t\le i)}\varkappa_S
\wedge dx_{s_1}\wedge\dots\wedge dx_{s_t}$ for some 
$\varkappa_S\in k(x_0)\otimes_k\Omega^{j-t}_{k/k_0}$. 

The $(1,\dots,1)$ weight subspace in $\bigoplus_{S=(0\le s_1
<\dots<s_t\le i)}k(x_0)\otimes_k\Omega^{j-t}_{k/k_0}\wedge 
dx_{s_1}\wedge\dots\wedge dx_{s_t}$ with respect to the action 
of homotheties $({\mathbb Q}^{\times})^{i+1}$ coincides with 
$x_0\cdot\Omega^{j-i}_{k/k_0}\wedge dx_1\wedge\dots\wedge dx_i
\oplus\Omega^{j-i-1}_{k/k_0}\wedge dx_0\wedge\dots\wedge dx_i$, 
so $\omega=\varkappa_1\wedge\eta+\varkappa_2\wedge d\eta$ 
for some $\varkappa_1\in\Omega^{j-i}_{k/k_0}$ and 
$\varkappa_2\in\Omega^{j-i-1}_{k/k_0}$, and therefore, 
$\varphi=\varkappa_1\wedge+\varkappa_2\wedge d$ (modulo 
$\Omega^1_{k_0}\wedge\Omega^{j-1}_F$). 

It follows from the above that ${\rm Hom}_G(k\otimes_{{\mathbb Z}}
\Omega^n_{F,{\rm exact}},\Omega^{<n}_{F/k_0})=0$. Check now that 
any $\varphi\in{\rm Hom}_G(\Omega^n_F,\Omega^{<n}_{F/k_0})$ factors 
through $\Omega^n_{F/k}\longrightarrow\Omega^{<n}_{F/k_0}$. For any 
$a\in k$ and $x_1,\dots,x_n\in F$ algebraically independent over $k$ 
one has $\varphi(a\cdot dx_1\wedge\dots\wedge dx_n)=\varphi(d(ax_1)
\wedge\dots\wedge dx_n)=0$, so $\varphi(x_1\cdot da\wedge dx_2
\wedge\dots\wedge dx_n)=0$. Clearly, $x_1\cdot da\wedge dx_2
\wedge\dots\wedge dx_n$ for all $a\in k$ are generators of the 
$G$-submodule $da\wedge\Omega^{n-1}_F$ of $\Omega^n_F$. This 
implies that ${\rm Hom}_G(\Omega^n_F,\Omega^{<n}_{F/k_0})=
{\rm Hom}_G(H^n_{{\rm dR}/k}(F),\Omega^{<n}_{F/k_0})$. 

More particularly, let $\varphi\in{\rm Hom}_G(\Omega^n_F,F)$. 
For any $L\subset F$ over which $F$ is of transcendence degree one 
${\rm Hom}_G(\Omega^n_F,F)\subseteq{\rm Hom}_{G_{F/\overline{L}}}
(\Omega^n_F,F)\subseteq{\rm Hom}_{G_{F/\overline{L}}}(\Omega^{n-1}_L
\otimes_L\Omega^1_F,F)$, so we may assume that $L=k$ and $n=1$. 
Let ${\rm Div}^{\circ}_{{\mathbb Q}}:=\lim\limits
_{_L\longrightarrow\phantom{_L}}{\rm Div}^{\circ}([L])_{{\mathbb Q}}$, 
${\rm Pic}^{\circ}_{{\mathbb Q}}:=\lim\limits
_{_L\longrightarrow\phantom{_L}}{\rm Pic}^{\circ}([L])_{{\mathbb Q}}$ 
and $H^1_{{\rm dR}/k,c}:=\lim\limits_{_L\longrightarrow\phantom{_L}}
H^1_{{\rm dR}/k}([L])_{{\mathbb Q}}$, where $[L]$ denotes a smooth proper 
model of $L$ over $k$, cf. \cite{repr}, p.182. Then 
${\rm Hom}_G({\rm Pic}^{\circ}_{{\mathbb Q}},F)={\rm Hom}_G
(H^1_{{\rm dR}/k,c},F)=0$, since if $\omega$ from ${\rm Pic}^{\circ}
_{{\mathbb Q}}$ or from $H^1_{{\rm dR}/k,c}$ is fixed by some $G_{F/L}$ 
and sent to $f\in L$ then ${\rm tr}_{L/L'}\omega\mapsto{\rm tr}
_{L/L'}f=[L:L']\cdot f$ for any $f\in L'\subseteq L$ purely transcendental 
over $k$, over which $L$ is algebraic, but ${\rm tr}_{L/L'}\omega=0$. 
So from the short exact sequences $0\longrightarrow 
H^1_{{\rm dR}/k,c}\longrightarrow H^1_{{\rm dR}/k}(F)\stackrel{Res}
{\longrightarrow}{\rm Div}^{\circ}_{{\mathbb Q}}\otimes k
\longrightarrow 0$ and $0\longrightarrow F^{\times}/
k^{\times}\stackrel{div}{\longrightarrow}{\rm Div}^{\circ}
_{{\mathbb Q}}\longrightarrow{\rm Pic}^{\circ}_{{\mathbb Q}}
\longrightarrow 0$, we get ${\rm Hom}_G(H^1_{{\rm dR}/k}(F),F)=0$. 

The form $\omega=\frac{dx_1}{x_1}\wedge\dots\wedge
\frac{dx_q}{x_q}$ is a generator of the $G$-module 
$\Omega^q_{F/k_0,\log}$. A $G$-homomorphism to 
$\Omega^{\bullet}_{F/k_0}$ sends it to an element $\omega'$ of 
$\left(\Omega^{\bullet}_{F/k_0}\right)^{{\rm Stab}_{\omega}}$. 
As ${\rm Stab}_{\omega}\supset U_{k(x_1,\dots,x_q)}$, 
$\left(\Omega^{\bullet}_{F/k_0}\right)^{{\rm Stab}_{\omega}}\subset
\Omega^{\bullet}_{k(x_1,\dots,x_q)/k_0}$. As $\omega'=\sum_I\eta_I
\wedge\frac{dx_I}{x_I}$ for some $\eta_I\in k(x_1,\dots,x_q)\otimes_k
\Omega^{\bullet}_{k/k_0}$ is fixed by $({\mathbb Q}^{\times})^q$ one 
has $\eta_I\in\Omega^{\bullet}_{k/k_0}$. As $\omega'$ is fixed by 
${\rm SL}_q({\mathbb Z})$ one has $\omega'=\eta_{1,\dots,q}\wedge
\frac{dx_1}{x_1}\wedge\dots\wedge\frac{dx_q}{x_q}$, where 
$\eta_{1,\dots,q}\in\Omega^{\bullet}_{k/k_0}$. \qed 

\vspace{4mm}

According to \cite{Hodge}, $H^q_{{\rm dR}/k}(F):={\rm coker}
[\Omega^{q-1}_{F/k}\stackrel{d}{\longrightarrow}\Omega^q
_{F/k,{\rm closed}}]=\lim\limits_{\longrightarrow}{\mathbb H}^q
(X,\Omega^{\bullet}_{X/k}(\log D))$, where $(X,D)$ runs over pairs 
consisting of a smooth proper variety $X$ with $k(X)\subset F$ and 
a normal crossing divisor $D$ on $X$ with smooth irreducible 
components. Moreover, the Hodge filtration on 
$\Omega^{\bullet}_{X/k}(\log D)$ induces a descending filtration 
on $H^q_{{\rm dR}/k}(F)$ by $k$-linear representations of $G_{F/k}$ 
with the graded quotients $\lim\limits_{\longrightarrow}H^p(X,
\Omega^{q-p}_{X/k}(\log D))$, where $(X,D)$ runs over pairs as above. 

The weight filtration on $\Omega^{\bullet}_{X/k}(\log D)$ induces 
an increasing filtration $W_{\bullet}H^q_{{\rm dR}/k}(F)$ on 
$H^q_{{\rm dR}/k}(F)$. In particular, $H^q_{{\rm dR}/k,c}(F):=
W_qH^q_{{\rm dR}/k}(F)$ is the image in $H^q_{{\rm dR}/k}(F)$ of 
$\lim\limits_{\longrightarrow}H^q_{{\rm dR}/k}(X)$, where $X$ 
runs over smooth proper models of subfields in $F$ of finite type 
over $k$. Clearly, this is an admissible representation over $k$. 
Again, the Hodge filtration on $\Omega^{\bullet}_{X/k}$ induces a 
descending filtration on $H^q_{{\rm dR}/k,c}(F)$ with the graded 
quotients $H^{p,q-p}_{F/k}=\lim\limits_{\longrightarrow}{\rm coker}
[H^{p-1}(D,\Omega^{q-p-1}_{D/k})\longrightarrow H^p(X,
\Omega^{q-p}_{X/k})]$, where $(X,D)$ runs over pairs as above. 
More particularly, $H^{q,0}_{F/k}=\Omega^q_{F/k,{\rm reg}}\subset 
H^q_{{\rm dR}/k,c}(F)$. 

\begin{proposition} Suppose that the cardinality of $k$ is at most 
continuum. Fix an embedding $\iota:k\hookrightarrow{\mathbb C}$ to 
the field of complex numbers. Then \begin{itemize} 
\item there is a non-canonical ${\mathbb Q}$-linear 
isomorphism $H^{p,q}_{F/k}\cong H^{q,p}_{F/k}$, and a 
${\mathbb C}$-anti-linear canonical isomorphism 
(depending on $\iota$) $H^{p,q}_{F/k}\otimes_{k,\iota}{\mathbb C}
\cong H^{q,p}_{F/k}\otimes_{k,\iota}{\mathbb C}$; 
\item the representation $H^n_{{\rm dR}/k,c}(F)$ is semi-simple 
for any $1\le n<\infty$. \end{itemize} \end{proposition}
{\it Proof.} \begin{itemize}\item The complexification of the 
projection $F^pH^{p+q}_{{\rm dR}/k}(X)\longrightarrow\hspace{-3mm}\to 
H^q(X,\Omega^p_{X/k})$ identifies the space $F^pH^{p+q}_{{\rm dR}/k}(X)
\otimes_{k,\iota}{\mathbb C}\cap\overline{F^qH^{p+q}_{{\rm dR}/k}(X)
\otimes_{k,\iota}{\mathbb C}}$ with $H^q(X,\Omega^p_{X/k})
\otimes_{k,\iota}{\mathbb C}$, where $F^{\bullet}$ is the Hodge 
filtration. Then the complex conjugation on $H^{p+q}(X_{\iota}
({\mathbb C}),{\mathbb C})=H^{p+q}(X_{\iota}({\mathbb C}),{\mathbb R})
\otimes_{{\mathbb R}}{\mathbb C}$ identifies $H^q(X,\Omega^p_{X/k})
\otimes_{k,\iota}{\mathbb C}$ with $H^p(X_{\iota}({\mathbb C}),
\Omega^q_{X_{\iota}({\mathbb C})})=H^p(X,\Omega^q_{X/k})
\otimes_{k,\iota}{\mathbb C}$. 
\item The semi-simplicity of $H^n_{{\rm dR}/k,c}(F)$ is equivalent to 
the semi-simplicity of the representation ${\mathbb C}\otimes_{k,\iota}
H^n_{{\rm dR}/k,c}(F)=\bigoplus_{p+q=n}{\mathbb C}\otimes_{k,\iota}
H^{p,q}_{F/k}$ of $G$. For the latter note that there is a positive 
definite $G$-equivariant hermitian form $({\mathbb C}\otimes_{k,\iota}
H^{p,q}_{F/k})\otimes_{id,{\mathbb C},c}({\mathbb C}\otimes_{k,\iota}
H^{p,q}_{F/k})\longrightarrow{\mathbb C}(\chi)$, where $c$ is 
the complex conjugation and $\chi$ is the modulus of $G$, given 
by $(\omega,\eta)=\int_{X_{\iota}({\mathbb C})}i^{n^2+2q}\omega
\wedge\overline{\eta}\cdot[G_{F/k(X)}]$ for any $\omega,\eta\in 
H^{p,q}_{{\rm prim}}(X_{\iota}({\mathbb C}))={\mathbb C}\otimes
_{k,\iota}H^q_{{\rm prim}}(X,\Omega^p_{X/k})\subset{\mathbb C}
\otimes_{k,\iota}H^{p,q}_{F/k}$. Here $H^{p,q}_{{\rm prim}}
(X_{\iota}({\mathbb C}))$ denotes the subspace orthogonal to the 
sum of the images of all Gysin maps $H^{p-1,q-1}(D)\longrightarrow 
H^{p,q}(X_{\iota}({\mathbb C}))$ for all desingularizations $D$ 
of all divisors on $X_{\iota}({\mathbb C})$. \qed \end{itemize} 

\vspace{4mm} 

Let us check that the sequence $0\longrightarrow H^1_{{\rm dR}/k}(X)
\oplus k\otimes(k(X)^{\times}/k^{\times})\stackrel{i}{\longrightarrow}
H^1_{{\rm dR}/k}(k(X))\stackrel{Res}{\longrightarrow}k\otimes
{\rm Pic}^0(X)\longrightarrow 0$ is exact, where for any divisor 
$x\in X^1$ the residue ${\rm res}_x\omega$ is defined, and gives rise 
to the map $H^1_{{\rm dR}/k}(k(X))\stackrel{Res}{\longrightarrow}
k\otimes{\rm Div}(X)$. 

The map $i$ sends a pair $(\omega,a\otimes f)$ to 
$\omega+ad\log f$, and has trivial kernel, since the residue map 
is trivial on $H_{{\rm dR}}^1(k(X)/k)$ and injective on 
$k\otimes(k(X)^{\times}/k^{\times})$.
If the residues of $\omega\in H^1_{{\rm dR}/k}(k(X))$ are zero 
then integration along a loop depends only on its homology 
class in $H_1(X,{\mathbb Q})$. There is an element $\eta$ of 
$H_{{\rm dR}/k}^1(X)$ with the same periods as $\omega$, 
so integration of $\omega-\eta$ along a path joining a fixed 
(rational) point with the variable one is independent of a 
chosen path, and defines a meromorphic (i.e. rational) function. 
This gives exactness in the middle term: ${\rm Im}i={\rm Ker}Res$. 

For any pair $D_1,D_2$ of algebraically equivalent effective 
divisors on $X$ there is a smooth projective curve $C$, and an 
effective divisor $D$ on $X\times C$, such that ${\rm pr}_X: 
D\longrightarrow X$ is generically finite and for some points 
$P,Q\in C$ one has $D_P-D_Q=D_1-D_2$. 

Since $\dim_k\Gamma(C,\Omega^1_C(P+Q))=\dim_k\Gamma(C,\Omega^1_C)+1$, 
$\Gamma(C,\Omega^1_C(P))=\Gamma(C,\Omega^1_C(Q))=\Gamma(C,\Omega^1_C)$, 
there exists a 1-form $\omega_{P,Q}$ such that $Res(\omega_{P,Q})=P-Q$. 

Set $\omega={\rm pr}_{X\ast}(({\rm pr}_C^{\ast}\omega_{P,Q})|_D)$. 
Then $Res(\omega)=D_1-D_2$. So the image of $Res$ contains the 
algebraically trivial part of the group of divisors with 
coefficients in $k$. This also shows that $\omega\in 
N_1\Omega^1_{{\rm dR}/k}(k(X))$. 
Since $Res$ commutes with restriction to a curve, $Res(\omega)\cdot
C=Res(\omega |_C)\in CH_0(X)$, $\deg(Res(\omega)\cdot C)=0$ by Cauchy 
theorem, $NS(X)_{{\mathbb Q}}\otimes CH_1(X)_{{\mathbb Q}}/hom
\longrightarrow{\mathbb Q}$ 
is non-degenerate, $Res(\omega)=0\in H^2(X,{\mathbb Q})$. 
Thus, the map  $Res$ is well-defined and surjective. \qed

\begin{proposition} \label{odin-formy} The representation $\Omega^1
_{F/k,{\rm closed}}$ admits the following description for any 
$1\le n\le\infty$. Let $H^1_{{\rm dR}/k,c}(F):=\ker[H^1_{{\rm dR}/k}
(F)\stackrel{{\rm Res}}{\longrightarrow}k\otimes{\rm Div}^{\circ}
_{{\mathbb Q}}]$; ${\rm Pic}^{\circ}_{{\mathbb Q}}:={\rm coker}
[F^{\times}/k^{\times}\stackrel{{\rm div}}{\longrightarrow}
{\rm Div}^{\circ}_{{\mathbb Q}}]$.\footnote{cf. p.182 before 
Proposition 3.11 of \cite{repr}} Then ${\rm Pic}^{\circ}_{{\mathbb Q}}
=\bigoplus_AA^{\vee}(k)\otimes_{{\rm End}(A)}(A(F)/A(k))$, where $A$ 
runs over the isogeny classes of simple abelian varieties over $k$. 
\begin{itemize} \item The maximal semi-simple subrepresentation of 
$G$ in $\Omega^1_{F/k,{\rm closed}}$ is canonically isomorphic to 
$$\bigoplus_A\Gamma(A,\Omega^1_{A/k})^{A(k)}\otimes_{{\rm End}(A)}
(A(F)/A(k))=(F/k)\oplus k\otimes(F^{\times}/k^{\times})
\oplus\Omega^1_{F/k,{\rm reg}},$$ where $A$ runs over the
isogeny classes of simple commutative algebraic $k$-groups. 
\item The maximal semi-simple subrepresentation of $G$ in 
$H^1_{{\rm dR}/k}(F)$ is canonically isomorphic to $$\bigoplus_A
H^1_{{\rm dR}/k}(A)\otimes_{{\rm End}(A)}(A(F)/A(k))=k\otimes
(F^{\times}/k^{\times})\oplus H^1_{{\rm dR}/k,c}(F),$$ where $A$ 
runs over the isogeny classes of simple commutative algebraic 
$k$-groups (with the zero summand corresponding to ${\mathbb G}_a$). 
\item The representation $H^1_{{\rm dR}/k}(F)/
(k\otimes(F^{\times}/k^{\times}))$ of $G$ is canonically 
isomorphic to $$\bigoplus_A[H^1_{{\rm dR}/k}(k(A))/
(k\otimes(k(A)^{\times}/k^{\times}))]\otimes_{{\rm End}(A)}(A(F)/A(k)),$$ 
where $A$ runs over the isogeny classes of simple abelian $k$-varieties. 
\end{itemize} \end{proposition} 
{\it Proof.} Follows from the above and from (evidently modified) 
Proposition 3.11 of \cite{repr}. \qed 

\begin{corollary} Let ${\mathbb Q}\subseteq k_0\subseteq k$. Then 
${\rm End}_{k_0[G]}(\Omega^1_{F/k,{\rm closed}})$ contains 
$${\rm Hom}_{k_0[G]}(\Omega^1_{F/k,{\rm closed}},\Omega^1_{F/k,{\rm reg}})\\ 
=\prod_A{\rm Hom}_{k_0}\left(H^1_{{\rm dR}/k}(k(A))/k\otimes(k(A)^{\times}/
k^{\times}),\Gamma(A,\Omega^1_{A/k})\right)$$ and ${\rm End}_{k_0[G]}
(H^1_{{\rm dR}/k}(F))$ contains (properly, if $k$ is transcendental over 
${\mathbb Q}$ and $k_0\subseteq\overline{{\mathbb Q}}$) $${\rm Hom}_{k_0[G]}
(H^1_{{\rm dR}/k}(F),H^1_{{\rm dR}/k,c}(F))\\ 
=\prod_A{\rm Hom}_{k_0}\left(H^1_{{\rm dR}/k}(k(A))/k\otimes
(k(A)^{\times}/k^{\times}),H^1_{{\rm dR}/k}(A)\right),$$ where $A$ runs 
over isogeny classes of simple abelian varieties over $k$. \end{corollary}
{\it Proof.} As there are no subobjects of $H^1_{{\rm dR}/k,c}(F)$ isomorphic 
to $F^{\times}/k^{\times}$, $${\rm Hom}_{k_0[G]}(H^1_{{\rm dR}/k}(F),
H^1_{{\rm dR}/k,c}(F))
={\rm Hom}_{k_0[G]}(H^1_{{\rm dR}/k}(F)/k\otimes(F^{\times}/k^{\times}),
H^1_{{\rm dR}/k,c}(F)),$$ so applying ${\rm Hom}_{k_0[G]}(H^1_{{\rm dR}/k}(F),-)$ to
$0\longrightarrow H^1_{{\rm dR}/k,c}(F)\longrightarrow H^1_{{\rm dR}/k}(F)
\stackrel{{\rm Res}}{\longrightarrow}k\otimes{\rm Div}^{\circ}_{{\mathbb Q}}
\longrightarrow 0$ we get the inclusion. On the other hand, the Gau\ss--Manin 
connection induces an embedding ${\rm Der}(k)\hookrightarrow
{\rm End}_{\overline{{\mathbb Q}}[G]}(H^1_{{\rm dR}/k}(F))$, which does 
not factor through ${\rm Hom}_G(H^1_{{\rm dR}/k}(F),H^1_{{\rm dR}/k,c}(F))$ 
unless it is zero, i.e. $k\subseteq\overline{{\mathbb Q}}$. \qed 

\vspace{4mm}

{\it Remark.} Clearly, the projection $\Omega^1_{F/k,{\rm closed}}
\longrightarrow\hspace{-3mm}\to\Omega^1_{F/k,{\rm closed}}/\Omega
^1_{F/k,{\rm reg}}$ is split. However, there is no natural splitting. 

\section{Alternative descriptions of ``homotopy invariant'' representations} 
\subsection{``Homotopy invariant'' representations as non-degenerate modules} 
\label{I-Hecke}
We set ${\bf D}_E:=\!\lim\limits_{\longleftarrow_U}\!E[G/U]$, where, 
for a field $E$ of characteristic zero, the inverse 
system is formed with respect to the projections 
$E[G/V]\stackrel{r_{VU}}{\longrightarrow}E[G/U]$ induced by inclusions 
$V\subset U$ of open subgroups in $H$.  
For any $\nu\in{\bf D}_E$, any $\sigma\in G$ and an open subgroup 
$U$ we set $$\nu(\sigma U):=\mbox{coefficient of $[\sigma U]$ of the 
image of $\nu$ in $E[G/U]$}.$$ The {\sl support} of $\nu$ is the 
minimal closed subset $S$ in $\widehat{H}$ such that $\nu(\sigma U)=0$ 
if $\sigma U\bigcap S=\emptyset$. Define a pairing ${\bf D}_E\times 
W\longrightarrow W$ for each smooth $E$-representation $W$ of $H$ by 
$(\nu,w)\longmapsto\sum_{\sigma\in G/V}\nu(\sigma V)\cdot\sigma w$, 
where $V$ is an arbitrary open subgroup in the stabilizer of $w$. When 
$W=E[G/U]$ this pairing is compatible with the projections $r_{VU}$, 
so we get a pairing ${\bf D}_E\times\lim\limits_{\longleftarrow_U}
\!E[G/U]\longrightarrow\lim\limits_{\longleftarrow_U}\!E[G/U]=
{\bf D}_E$, and thus an associative multiplication ${\bf D}_E
\times{\bf D}_E\stackrel{\ast}{\longrightarrow}{\bf D}_E$ 
extending the convolution of compactly supported measures. 

If $n=\infty$ then the action of the associative algebra ${\bf D}_E$ 
on any object of ${\mathcal I}_G(E)$ factors through the action of 
its quotient ${\bf C}_E:=\lim\limits_{\longleftarrow_L}C_L$, since 
the morphism $E[G/G_{F/L}]\otimes_EW^{G_{F/L}}\longrightarrow W$ of 
representations of $G$ factors through ${\mathcal I}(E[G/G_{F/L}]
\otimes_EW^{G_{F/L}})=C_L\otimes W^{G_{F/L}}\longrightarrow W$. 

For any compact subgroup $U$ in $G$ the action of its Hecke algebra 
${\mathcal H}_E(U):=h_U\ast{\bf D}_E\ast h_U$ on $W^U$ factors through 
the action of its quotient ${\mathcal C}_E(U):=h_U\ast{\bf C}_E\ast h_U$ 
in ${\bf C}_E$ for any $W\in{\mathcal I}_G(E)$. E.g., if $F^U$ is purely 
transcendental over $L$ and $L$ is of finite type $k$ then 
${\mathcal C}_E(U)=C_L^U\otimes E$. 

Let $\mathcal{H}_{{\mathcal I}}:=
\lim\limits_{_K\longrightarrow\phantom{_K}}\lim\limits
_{\phantom{_L}\longleftarrow_L}C_L^{G_{F/K}}$ be the associative 
idempotented algebra without unity. The images $h_K$ of the Haar 
measures on $G_{F/K}$ for purely transcendental extensions $K$ of 
subfields of finite type over $k$ in $F$ over which $F$ is algebraic, 
are projectors in the algebra $\mathcal{H}_{{\mathcal I}}$. Then 
the category ${\mathcal I}_G$ is equivalent to the category of 
non-degenerate modules over $\mathcal{H}_{{\mathcal I}}$, i.e. such 
modules $W$ that $W=\mathcal{H}_{{\mathcal I}}W$. The algebra 
$\mathcal{H}_{{\mathcal I}}$ is isomorphic to the Hecke algebra 
(of locally invariant measures with compact support) of neither 
locally compact group, since any, e.g. finite-dimensional, subspace 
in $\lim\limits_{\phantom{_L}\longleftarrow_L}C_L^G$ is a left ideal 
in $\mathcal{H}_{{\mathcal I}}$, which never happens in the Hecke 
algebras.\footnote{Let $H$ be a locally compact group. If there is 
a non-zero finite-dimensional left ideal ${\mathfrak a}$ in the 
Hecke algebra of $H$ then the common support of the measures in 
${\mathfrak a}$ is compact and left-invariant, and therefore, $H$ 
is compact. Then the smooth representations of $H$ are semi-simple.} 

\subsection{A locally compact ``dense subgroup'' of $G$ 
and a description of ${\mathcal I}_G$, etc.} 
In the case $n=\infty$ the category ${\mathcal I}_G$ admits also a 
description in terms of a locally compact group. 

For a descending sequence $L_{\bullet}=(L_1\supset L_2\supset 
L_3\supset\dots)$ of subfields in $F$ set $H=H_{L_{\bullet}}=
\bigcup\limits_{m\ge 1}G_{F/L_m}$. We take the subgroups 
$G_{F/L_1(S)}\subset H$ for all finite subsets $S$ in $F$ 
as a base of open subgroups. 

\begin{itemize}\item We want $H$ to be a dense ``subgroup'' of $G$, so 
we ask that $\bigcap_{m\ge 1}\overline{L_m}=k$. This implies that the
forgetful functor ${\mathcal S}m_G\longrightarrow H\mbox{-}{\rm mod}$
is fully faithful.\footnote{{\it Proof.} Let $W,W'\in{\mathcal S}m_G$, 
$\alpha\in{\rm Hom}_H(W,W')$, $v\in W$ and $\sigma\in G$. Let $U$ be 
the common stabilizer of $v$ and $\alpha(v)$. Choose some $\sigma'\in 
H\bigcap\sigma U$. Then $\alpha(\sigma v)=\alpha(\sigma'v)=
\sigma'\alpha(v)=\sigma\alpha(v)$. \qed} 
\item Further, we want $H$ to be locally compact and therefore we 
ask $F$ to be of finite transcendence degree over $L_1$. Then $H$ 
is indeed locally compact, but not unimodular. 
\item If $L'_{\bullet}$ is an infinite subset in $L_{\bullet}$ then 
clearly $H_{L_{\bullet}}=H_{L'_{\bullet}}$. We want the topologies 
on $H_{L_{\bullet}}$ and on $H_{L'_{\bullet}}$ to be the same. For 
that we ask $L_1$ to be of finite type over $L_m$ for all $m>1$. 
\item The inclusion of $H$ into $G$ is a continuous homomorphism, so the 
forgetful functor ${\mathcal S}m_G\longrightarrow H\mbox{-}{\rm mod}$ 
factors through ${\mathcal S}m_G\longrightarrow{\mathcal S}m_H$.
It admits a right adjoint $W\mapsto\lim\limits_{\longrightarrow}
\bigcap_{m\ge 1}W^{G_{F/LL_m}}$, where $L$ runs over subfields of 
finite type over $k$ in $F$. The $G$-action is defined as follows. 
If $w\in W^{G_{F/LL_m}}$ and $\sigma\in G$ then $\sigma w:=\sigma'w$, 
where $\sigma'\in H$ and $\sigma'|_L=\sigma|_L$. Clearly, this is 
independent of $\sigma'$. 
\item Suppose that $L_j$ is purely transcendental over $L_{j+1}$
for any $j\gg 1$. As any admissible representation of $G$ 
is `homotopy invariant', the forgetful functor induces 
${\mathcal A}dm_G\longrightarrow{\mathcal A}dm_H$. \end{itemize}

In particular, the effective motives modulo numerical equivalence 
form a full subcategory in the category of graded semi-simple 
admissible $H$-modules. Note, that as the category of graded 
semi-simple admissible $H$-modules of finite length is self-dual, 
arbitrary motives (not necessarily effective) modulo numerical 
equivalence can probably be realized in that category. 

Let ${\mathcal I}_H$ be the full subcategory in ${\mathcal S}m_H$ 
whose objects $W$ satisfy the ``homotopy invariance'' condition 
$W^{G_{F/LL_m}}=W^{G_{F/LL_m(S)}}$ for any $m\ge 1$, any extension 
$L$ of $k$ of finite type in $F$ and any transcendence basis $S$ 
of $F$ over $LL_m$. 

{\sc Example.} Choose a transcendence basis $\{x_1,x_2,\dots\}$ of 
$F$ over $k$ and set $L_m=k(x_m,x_{m+1},\dots)$. 

Geometrically, this corresponds to an inverse system of 
infinite-dimensional irreducible varieties given by 
finite systems of equations. They are related by dominant 
morphisms affecting only finitely many coordinates. 

It follows from Lemma 2.15 of \cite{repr} that 
$H=G_{F/\overline{L_2}}\cdot H^{\circ}$. 

\begin{proposition} The forgetful functor ${\mathcal S}m_G
\longrightarrow{\mathcal S}m_H$ induces the following equivalences 
of categories: ${\mathcal I}_G\stackrel{\sim}{\longrightarrow}
{\mathcal I}_H$ and ${\mathcal A}dm_G\stackrel{\sim}
{\longrightarrow}{\mathcal I}_H\cap{\mathcal A}dm_H$. \end{proposition}
{\it Proof.} We need to construct an inverse functor 
${\mathcal I}_H\longrightarrow{\mathcal I}_G$. In particular, 
for a given $W\in{\mathcal I}_H$, $v\in W$ and $\sigma\in G$, 
we want to define $\sigma v$. 

There exist a subfield $L\subset F$ of finite type over $k$ and an 
integer $m\ge 1$ such that the stabilizer of $v$ contains $G_{F/LL_m}$. 
Let $LL_m=L'L_{m'}$, where $L'\subset F$ is of finite type 
over $k$, and $L'$ and $L_{m'}$ are algebraically independent over $k$. 

Let $N>m'$ be an integer such that $L'\sigma(L')$ and $L_N$ are 
algebraically independent over $k$. Take any $\sigma'\in G_{F/L_N}$ 
such that $\sigma'|_{L'}=\sigma|_{L'}$ and set $\sigma v:=\sigma'v$. 
One has $v\in W^{G_{F/L'L_{m'}}}=W^{G_{F/L'L_N}}$, so $\sigma v$ 
is independent of particular choices of $N$ and of $\sigma'$. 

Now we check independence of $L'$. Suppose that $v\in W^{G_{F/L'L_{m'}}}
\cap W^{G_{F/L''L_{m''}}}$. Since $v\in W^{G_{F/L'L''L_{m'+m''}}}$, it 
suffices to treat the case $L'\subseteq L''$. As above, we choose an 
integer $N>m''$ such that $L''\sigma(L'')$ and $L_N$ are algebraically 
independent over $k$, and some  $\sigma''\in G_{F/L_N}$ such that 
$\sigma''|_{L''}=\sigma|_{L''}$. Then $\sigma''$ can also 
serve as a $\sigma'$, i.e., $\sigma''v=\sigma'v$. 

This gives us a map $G\times W\longrightarrow W$. Clearly, 
this is a linear action, and the stabilizer of $v$ contains 
the open subgroup $G_{F/L'}$, and thus, $W$ becomes 
an object of ${\mathcal I}_G$. \qed 

{\it Remarks.} 1. There are admissible representations of $H$ outside 
of ${\mathcal I}_H$, e.g. ${\mathbb Q}(\rho)$ for any non-trivial 
character $\rho$ of $H$. 

2. ${\mathcal I}_H$ is closed under subquotients and direct 
products,\footnote{the direct product of a family of smooth 
representations is the smooth part of their set-theoretic 
direct product} but not under extensions in ${\mathcal S}m_H$. 
As any morphism from $W\in{\mathcal S}m_H$ to an object of 
${\mathcal I}_H$ factors through the canonical map to the 
direct product over all morphisms from $W$ to representatives 
of all isomorphism classes in ${\mathcal I}_H$, there is 
a functor ${\mathcal I}:{\mathcal S}m_H\longrightarrow
{\mathcal I}_H$ left adjoint to the inclusion 
${\mathcal I}_H\hookrightarrow{\mathcal S}m_H$. 

3. ${\mathcal A}dm_H$ is a Serre subcategory in ${\mathcal S}m_H$.

\end{document}